\documentclass[11pt]{article}

\usepackage{latexsym, amsmath, amssymb, eufrak} 

\usepackage{mathptmx} 
\usepackage{avant} 
\usepackage{authblk}
\usepackage{multirow}

\input amssym.def
\input amssym

\usepackage{graphicx}
\usepackage{caption}
\usepackage{float}

\usepackage[a4paper, margin=1in]{geometry}

\def\dj{d\kern-0.3em\char"16\kern-0.08em}
\def\Dj{\mbox{\raise0.3ex\hbox{-}\kern-0.3em D}}

\usepackage[bottom]{footmisc}

\usepackage{theorem}

\newtheorem{theorem}{Theorem}[section]

\theoremstyle{definition}
\newtheorem{definition}[theorem]{Definition}
\theoremstyle{remark}
\newtheorem{remark}[theorem]{Remark}

\input amssym.def
\input amssym

\DeclareMathOperator{\diag}{diag}
\DeclareMathOperator{\Elo}{Elo}

\def\R{\mathbb R}

\def\be{\begin{equation}}
\def\ee{\end{equation}}

\providecommand{\keywords}[1]{\textbf{\textit{Keywords:}} #1}

\providecommand{\subjclass}[1]{\textbf{\textit{Mathematics Subject Classification:}} #1}

\begin{document}

\title{Bayesian statistics approach\\ to chess engines optimization}

\author[1]{Ivan Ivec}
\author[2]{Ivana Vojnovi\' c}

\affil[1]{Faculty of Metallurgy, University of Zagreb, Aleja narodnih heroja 3, Sisak, Croatia\\

iivec@simet.unizg.hr}

\affil[2]{Department of Mathematics and Informatics, Faculty of Sciences, University of Novi Sad, Serbia\\

ivana.vojnovic@dmi.uns.ac.rs}

\maketitle

\begin{abstract}
We develop a new method for stochastic optimization using the Bayesian statistics approach. More precisely, we optimize parameters of chess engines as those data are available to us, but
the method should apply to all situations where we want to optimize a certain gain/loss function which has no analytical form and thus cannot be measured directly but only by comparison of two parameter sets. We also
experimentally compare the new method with the famous SPSA method.
\end{abstract}



\let\thefootnote\relax\footnote{This work has been supported in part by the Croatian Science Foundation
under the project 2449 MiTPDE. The second author acknowledges the financial support of the Ministry of Education, Science and Technological Development of the Republic of Serbia (Grant No. 451-03-68/2020-14/ 200125).}
\let\thefootnote\relax\footnote{We are thankful to Professor Nata\v sa Krklec - Jerinki\' c  for her feedback and valuable discussions.}

\section{Introduction}

The simultaneous perturbation stochastic approximation (SPSA) method for multivariate optimization problems was introduced in \cite{spa}, and it found numerous applications in engineering and the physical and social sciences.
The method was more fully analyzed in \cite{spb}, while in \cite{spc}, it was given a simple step-by-step guide to the implementation of SPSA, with references to various applications.

More recently, the SPSA method was implemented in \cite{spsa} as a tuner of parameters of the Stockfish chess engine. Stockfish is an open source project \cite{sf}, and it is currently the strongest
chess engine in the world (winning the last five seasons of the famous TCEC competition \cite{tcec} and being on the top of the CCRL rating lists \cite{ccrl}), much stronger than the engine used by the famous Deep Blue computer, which defeated a reigning chess world champion, Garry Kasparov, in 1997.
Over the years, the SPSA method brought many successful patches to Stockfish and other chess engines.

However, tuning parameters of chess engines is computationally very demanding, and many tuning attempts are unsuccessful. In this paper, we try to use the Bayesian statistics approach to find a method that will be better, at least in some situations. We choose this approach because Bayesian inference is known for giving good posterior distributions from pretty bad prior distributions if the prior's support and shape are reasonably selected. For more details and references regarding Bayesian optimization, see \cite{fra}.

In the next section, we provide some theoretical justifications for our method, which we call the BSPSA method because of the Bayesian approach and because of a similarity with the SPSA method.
In fact, a difference from SPSA method will be in a way of updating parameters' approximations. We want to optimize parameters $\theta=(\theta^{(1)},\theta^{(2)},\ldots,\theta^{(n)})\in\R^n$, and because the game of chess is a stochastic process, that actually means that we are looking for a distribution of the optimal value of $\theta$. In other words, we try to obtain $\theta$ that will work best in most situations, and the mean value of the distribution will be a candidate for that. SPSA updates parameters by formula
\[\theta_{k+1}^{(i)}=\theta_k^{(i)}+\frac{a_k^{(i)}}{\Delta_k^{(i)}c_k^{(i)}}\cdot w_k\,, \]
where
$a_k^{(i)},\Delta_k^{(i)}$ and $c_k^{(i)}$ are appropriate sequences and $w_k$ is the result of an experiment (two-game match in our case) -- see the beginning of the third section for more details. BSPSA will update parameters using the famous Bayes' theorem for the conditional probability, i.e. the formula (\ref{bay}) given below.
We provide experimental results in the third section, comparing SPSA and BSPSA.
Although we stick to tuning the Stockfish parameters, we believe the method could also be applied to other situations where a gain/loss function is not measured directly but only by comparing two parameter sets. We discuss such possibilities in the concluding section.

\section{Theoretical development of a new method}

For simplicity, we first look at a single parameter $\theta\in\R$ that we want to optimize by playing many chess games.
The central limit theorem guarantees that the distribution of experimentally obtained optimal values using $n$ games (for $n$ large)
will be approximately normal. Of course, it is always possible that we encounter an ill-behaved parameter whose optimal value will have an infinite mean or infinite variance and where the central limit theorem does not apply. Still, such parameters are not useful in the code of our engine, and we strongly try to avoid them. Therefore, we construct a sequence $(\theta_k)$, hopefully converging to the optimal value $\theta_0$,
by starting with $N(\theta_k,\,s_k^2)$ normal prior distribution and calculating $N(\theta_{k+1},\,s_{k+1}^2)$ posterior distribution.

So, assuming that we have already obtained $\theta_k$, we play two games between an engine $E_1$ that uses $\theta=\theta_k + c_k$
and an engine $E_2$ that uses $\theta=\theta_k - c_k$, where $c_k>0$ will be chosen in the same way as in the SPSA method or in some other way,
whichever proves to be more efficient. The result of this two-game match, from the perspective of the engine $E_1$, we denote by $w_k$ (with possible
values $2,1,0,-1,-2$, coming from one point for a win, zero points for a draw and $-1$ for a loss).

Now, if $f(w_k|\theta)$ is a conditional observation distribution of $w_k$, assuming that $\theta$ is given, then the Bayesian inference starting from the prior distribution with the density
\[\pi(\theta)\sim e^{-\frac12\big(\frac{\theta-\theta_k}{s_k}\big)^2} \]
gives the posterior distribution with the density
\begin{equation}\label{bay}
\pi(\theta|w_k)\sim f(w_k|\theta)\pi(\theta)\,, 
\end{equation}
where by $a\sim b$ we denote that $\frac ab$ is a constant.

In what follows, we show that, under reasonable assumptions, we can approximate $f(w_k|\theta)$ with the density of the normal distribution with the mean of the form $A\theta+B$,
for some constants $A>0$, $B\in\R$. Although $w_k$ can obtain only $5$ values, approximation with normal distribution is viable because of the known fact from the Bayesian inference that it makes no difference whether we analyze observations one at a time in sequence using the posterior after the previous step as the prior for the next step, or whether we analyze all observations together in a single step starting with our initial prior \cite[Sections 6.1, 11.1]{bol}. More precisely, that is why we decide to update our belief about the parameter after every two games instead of playing more games with the same parameter.

In that way, the posterior distribution $\pi(\theta|w_k)$ will also be normal (at least approximately), and we obtain explicit formulae for its calculation,
thus avoiding the numerical integration needed to get a constant in (\ref{bay}).

By $WR$ we denote a win rate, i.e.~a ratio between wins/points obtained by the engine $E_1$ and wins/points obtained by the engine $E_2$. Now we take an assumption that the strength of an engine (expressed as the number of wins/points obtained)
is proportional to the density of the normal distribution with mean $\theta$ (as we assume that $\theta$ is optimal) and standard deviation $\sigma$ (there is no clear justification for this assumption, but ideally we would like to construct an engine which parameters are independent and identically distributed with a finite variance, in which case the central limit theorem would be a kind of justification again), thus obtaining
\[WR= e^{-\frac12\big(\frac{\theta_k+c_k-\theta}{\sigma}\big)^2+\frac12\big(\frac{\theta_k-c_k-\theta}{\sigma}\big)^2}=e^{\frac{2c_k(\theta-\theta_k)}{\sigma^2}}\,. \]
The win percentage is then 
\[WP=\frac{WR}{WR+1}\,, \]
and thus we get
\begin{equation}\label{wkfa}
w_k\approx 2WP-2(1-WP)=2\cdot\frac{WR-1}{WR+1}=2\cdot\frac{e^{\frac{2c_k(\theta-\theta_k)}{\sigma^2}}-1}{e^{\frac{2c_k(\theta-\theta_k)}{\sigma^2}}+1}\,. 
\end{equation}
Here we approximate $w_k$ with its expectation ($2WP$ comes from two games with the win percentage $WP$ and one point for each win, while $-2(1-WP)$ comes from two games with the loss percentage $1-WP$ and $-1$ point for each loss).

As the most difficult task in practice is to optimize a parameter that is already pretty close to the optimum, we can take a linear approximation of the right-hand side in (\ref{wkfa}) around $\theta=\theta_k$ to obtain
\begin{equation}\label{wkapprox}
w_k\approx\frac{2c_k(\theta-\theta_k)}{\sigma^2}\,.
\end{equation}
Finally, we use (\ref{wkapprox}) to assume that $f(w_k|\theta)$ is approximately equal to the density of the normal distribution with mean $\frac{2c_k(\theta-\theta_k)}{\sigma^2}$
and standard deviation $\tau$. We also assume that standard deviations $\sigma$ and $\tau$ are constants not depending on $\theta$ because we are satisfied with high precision only around $\theta=\theta_k$.

For the brevity of calculations we denote $A=\frac{2c_k}{\sigma^2}$, $B=\frac{-2c_k\theta_k}{\sigma^2}$ and from (\ref{bay}) we obtain
\begin{align}
\pi(\theta|w_k)&\sim e^{-\frac12\big(\frac{w_k-A\theta-B}{\tau}\big)^2}\cdot e^{-\frac12\big(\frac{\theta-\theta_k}{s_k}\big)^2} \nonumber \\
&=e^{-\frac12\Big(\frac{w_k^2+A^2\theta^2+B^2-2A\theta w_k-2Bw_k+2AB\theta}{\tau^2}+\frac{\theta^2-2\theta_k\theta+\theta_k^2}{s_k^2}\Big)} \nonumber \\
&=e^{-\frac12\cdot\frac{(A^2s_k^2+\tau^2)\theta^2-2(Aw_k s_k^2-ABs_k^2+\theta_k\tau^2)\theta+\ldots}{s_k^2\tau^2}} \nonumber \\
&\sim e^{-\frac12\cdot\frac{A^2s_k^2+\tau^2}{s_k^2\tau^2}\Big(\theta-\frac{Aw_k s_k^2-ABs_k^2+\theta_k\tau^2}{A^2s_k^2+\tau^2}\Big)^2}\,. \nonumber
\end{align}
We have thus obtained that the posterior distribution is also a normal distribution with the mean
\[\theta_{k+1}=\frac{Aw_k s_k^2-ABs_k^2+\theta_k\tau^2}{A^2s_k^2+\tau^2} \]
and the variance
\[s_{k+1}^2=\frac{s_k^2\tau^2}{A^2s_k^2+\tau^2}\,. \]
After inserting $A=\frac{2c_k}{\sigma^2}$, $B=\frac{-2c_k\theta_k}{\sigma^2}$ we obtain final formulae which
we use to update our belief about the parameter after each two-game match:
\begin{align}
\theta_{k+1}&=\theta_k+\frac{2c_k s_k^2\sigma^2}{4c_k^2 s_k^2+\tau^2\sigma^4}\cdot w_k\,, \label{tup} \\
s_{k+1}^2&=\frac{s_k^2\tau^2\sigma^4}{4c_k^2 s_k^2+\tau^2\sigma^4}\,. \label{sup}
\end{align}

\begin{remark}\label{rma}
In the limit $\tau\to 0$ the formula (\ref{tup}) reads
\[\theta_{k+1}=\theta_k+\frac{\sigma^2}{2c_k}\cdot w_k\,, \]
which is exactly the formula used in the SPSA method, with the difference that SPSA uses $a_k=\frac{a}{(A+k)^\alpha}$ instead of $\frac{\sigma^2}2$.
This means that uncertainty of the two-game match,
contained in the parameter $\tau$, is our replacement for the decreasing sequence $a_k$ used in the SPSA method. Also, it is obvious from the formula (\ref{sup}) that $0<s_{k+1}<s_k$, and thus $\lim\limits_{k\to\infty}s_k$ exists. By passing to the limit in (\ref{sup}) we get $\lim\limits_{k\to\infty}s_k=0$ if $\lim\limits_{k\to\infty}c_k\ne 0$, which is a good sign to expect the convergence of $\theta_k$ to the optimal value. At first sight, it seems that $c_k$ used in the SPSA method (and so also in our method) does not satisfy this condition. However, we first choose the number of iterations $N$ and then we can choose $c_N$ to be big enough and we have $c_k\geq c_N>0$ for all $k\in\{1,\ldots,N\}$. In addition, this is another reason (besides those mentioned in the concluding section) for looking for a better choice of $c_k$ in future work.
\end{remark}

\begin{remark}
If we want to optimize more parameters $\theta=(\theta^{(1)},\theta^{(2)},\ldots,\theta^{(n)})$, then assuming that we have already obtained $\theta_k$, we play two games between an engine $E_1$ that uses $\theta=\theta_k + \Delta_k c_k$
and an engine $E_2$ that uses $\theta=\theta_k - \Delta_k c_k$, where each component of $c_k$ is positive and each component of $\Delta_k$ is using a Bernoulli $\pm 1$ distribution with probability $\frac12$ for each $\pm 1$ outcome. This choice satisfies conditions introduced in \cite{spb} and has already been used in \cite{spsa}. In our prior distribution we assume that all parameters are independent:
\[\pi(\theta)\sim e^{-\frac12\sum\limits_{i=1}^n\Big(\frac{\theta^{(i)}-\theta_k^{(i)}}{s_k^{(i)}}\Big)^2}\,, \]
and similarly as above we obtain
\[WR = e^{\sum\limits_{i=1}^n\frac{2\Delta_k^{(i)}c_k^{(i)}(\theta^{(i)}-\theta_k^{(i)})}{[\sigma^{(i)}]^2}} \]
and
\[w_k\approx\sum\limits_{i=1}^n\frac{2\Delta_k^{(i)}c_k^{(i)}(\theta^{(i)}-\theta_k^{(i)})}{[\sigma^{(i)}]^2}\,. \]
This means that we can assume that $f(w_k|\theta)$ is approximately equal to the density of the normal distribution with mean $\sum\limits_{i=1}^n\frac{2\Delta_k^{(i)}c_k^{(i)}(\theta^{(i)}-\theta_k^{(i)})}{[\sigma^{(i)}]^2}$
and standard deviation $\tau$. Using notation
\begin{equation}\label{not}
A^{(i)}=\frac{2\Delta_k^{(i)}c_k^{(i)}}{[\sigma^{(i)}]^2},\quad B^{(i)}=\frac{-2\Delta_k^{(i)}c_k^{(i)}\theta_k^{(i)}}{[\sigma^{(i)}]^2}
\end{equation}
we finally obtain the posterior distribution
\begin{equation}\label{multipost}
\pi(\theta|w_k)\sim e^{-\frac12\Big(\frac{w_k-\sum\limits_{i=1}^n (A^{(i)}\theta^{(i)}+B^{(i)})}{\tau}\Big)^2}\cdot e^{-\frac12\sum\limits_{i=1}^n\Big(\frac{\theta^{(i)}-\theta_k^{(i)}}{s_k^{(i)}}\Big)^2}\,.
\end{equation}
We see that this is again a multivariate normal distribution but not with independent variables. In what follows, we first try to approximate it with an independent case and then perform 
a full calculation using a covariance matrix. 
\end{remark}

Formula (\ref{multipost}) could be rewritten as
\[\pi(\theta|w_k)\sim e^{-\frac12\sum\limits_{i=1}^n\big([\theta^{(i)}]^2\Big(\frac{[A^{(i)}]^2}{\tau^2}+\frac1{[s_k^{(i)}]^2}\Big) - 2\theta^{(i)}\Big(\frac{A^{(i)}w_k-A^{(i)}\sum\limits_{j=1}^nB^{(j)}-\frac12 A^{(i)}\sum\limits_{j\neq i} 
A^{(j)}\theta^{(j)}}{\tau^2}+\frac{\theta_k^{(i)}}{[s_k^{(i)}]^2}\Big)\big)}\,, \]
and if we use an approximation 
\begin{equation}\label{indapprox}
\sum\limits_{j\neq i} A^{(j)}\theta^{(j)}\approx\sum\limits_{j\neq i} A^{(j)}\theta_k^{(j)}=-\sum\limits_{j\neq i} B^{(j)}\,,
\end{equation} 
we obtain
\begin{align}
\pi(\theta|w_k)&\sim e^{-\frac12\sum\limits_{i=1}^n\big([\theta^{(i)}]^2\Big(\frac{[A^{(i)}]^2}{\tau^2}+\frac1{[s_k^{(i)}]^2}\Big) - 2\theta^{(i)}\Big(\frac{A^{(i)}w_k-A^{(i)}B^{(i)}-\frac12 A^{(i)}\sum\limits_{j\neq i}B^{(j)}}{\tau^2}+\frac{\theta_k^{(i)}}{[s_k^{(i)}]^2}\Big)\big)} \nonumber \\
&\sim e^{-\frac12\sum\limits_{i=1}^n\frac{[A^{(i)}]^2[s_k^{(i)}]^2+\tau^2}{[s_k^{(i)}]^2\tau^2}\Big(\theta^{(i)}-\frac{[s_k^{(i)}]^2A^{(i)}\big(w_k-B^{(i)}-\frac12 \sum\limits_{j\neq i}B^{(j)}\big)+\theta_k^{(i)}\tau^2}{[A^{(i)}]^2[s_k^{(i)}]^2+\tau^2}\Big)^2}\,. \nonumber
\end{align}
We have thus obtained (for any $i=1,2,\ldots ,n$)
\[\theta_{k+1}^{(i)} = \frac{[s_k^{(i)}]^2A^{(i)}\big(w_k-B^{(i)}-\frac12 \sum\limits_{j\neq i}B^{(j)}\big)+\theta_k^{(i)}\tau^2}{[A^{(i)}]^2[s_k^{(i)}]^2+\tau^2} \]
and
\[[s_{k+1}^{(i)}]^2=\frac{[s_k^{(i)}]^2\tau^2}{[A^{(i)}]^2[s_k^{(i)}]^2+\tau^2}\,. \]
Finally, using (\ref{not}) we get
\begin{align}
\theta_{k+1}^{(i)}&=\theta_k^{(i)}+\sum\limits_{j\neq i}\frac{2\Delta_k^{(i)}\Delta_k^{(j)}c_k^{(i)}c_k^{(j)}[s_k^{(i)}]^2[\sigma^{(i)}]^2}{[\sigma^{(j)}]^2(4[c_k^{(i)}]^2[s_k^{(i)}]^2+\tau^2[\sigma^{(i)}]^4)}\cdot \theta_k^{(j)}+\frac{2\Delta_k^{(i)}c_k^{(i)}[s_k^{(i)}]^2[\sigma^{(i)}]^2}{4[c_k^{(i)}]^2[s_k^{(i)}]^2+\tau^2[\sigma^{(i)}]^4}\cdot w_k\,, \label{tindup} \\
[s_{k+1}^{(i)}]^2&=\frac{[s_k^{(i)}]^2\tau^2[\sigma^{(i)}]^4}{4[c_k^{(i)}]^2[s_k^{(i)}]^2+\tau^2[\sigma^{(i)}]^4}\,. \label{sindup}
\end{align}

Moreover, we expect better results if we do not use the approximation (\ref{indapprox}). To do so, we take a general multivariate normal prior
\[\pi(\theta)\sim e^{-\frac12(\theta - \theta_k)^\top S_k^{-1}(\theta - \theta_k)}\,, \]
with a positive definite covariance matrix $S_k$ to rewrite (\ref{multipost}) in the form
\[\pi(\theta|w_k)\sim e^{-\frac12\Big(\frac{w_k-\sum\limits_{i=1}^n (A^{(i)}\theta^{(i)}+B^{(i)})}{\tau}\Big)^2}\cdot e^{-\frac12(\theta - \theta_k)^\top S_k^{-1}(\theta - \theta_k)}\,. \]
To simplify we use a substitution $\tilde{\theta}=\theta-\theta_k$ to obtain
\[\pi(\tilde\theta)\sim e^{-\frac12\tilde\theta^\top S_k^{-1}\tilde\theta} \]
and
\begin{equation}\label{updA}
\pi(\tilde\theta|w_k)\sim e^{-\frac12\Big(\frac{w_k-\sum\limits_{i=1}^n A^{(i)}\tilde\theta^{(i)}}{\tau}\Big)^2}\cdot e^{-\frac12\tilde\theta^\top S_k^{-1}\tilde\theta}\,. 
\end{equation}
We want to show that $\pi(\tilde\theta|w_k)$ could be written in the form
\begin{equation}\label{updB}
\pi(\tilde\theta|w_k)\sim e^{-\frac12(\tilde\theta-b_k)^\top S_{k+1}^{-1}(\tilde\theta-b_k)}\,, 
\end{equation}
for some vector $b_k$ and some matrix $S_{k+1}^{-1}$. We obtain that by equating the coefficients of polynomials in the exponents of (\ref{updA}) and (\ref{updB}). We denote elements of $S_k^{-1}$ as $(s_{ij})$ and elements of $S_{k+1}^{-1}$ as $(t_{ij})$, and after equating the coefficients in front of $[\tilde\theta^{(i)}]^2$ we get
\begin{equation}\label{diag}
t_{ii}=s_{ii}+\frac{[A^{(i)}]^2}{\tau^2}\,,
\end{equation}
while after equating the coefficients in front of $\tilde\theta^{(i)}\tilde\theta^{(j)}$ (for $i\ne j$) we get
\begin{equation}\label{nondiag}
t_{ij}=s_{ij}+\frac{A^{(i)}A^{(j)}}{\tau^2}\,.
\end{equation}
So, we obtain simple rules for updating elements of $S_k^{-1}$ given by (\ref{diag}) and (\ref{nondiag}). Actually, (\ref{nondiag}) covers also (\ref{diag}) by taking $i=j$.

Obtaining $b_k$ is not straightforward. Namely, by equating the coefficients in front of $\tilde\theta^{(i)}$ we get
\[t_{ii}b_k^{(i)}+\sum_{j\ne i}t_{ij}b_k^{(j)}=\frac{A^{(i)}w_k}{\tau^2}\,,\]
which is a system of linear equations, which written in the matrix form reads
\begin{equation}\label{parupd}
S_{k+1}^{-1}\cdot\begin{bmatrix}b_k^{(1)} \\ b_k^{(2)} \\ \vdots \\ b_k^{(n)}\end{bmatrix}=\frac{w_k}{\tau^2}\begin{bmatrix}A^{(1)} \\ A^{(2)} \\ \vdots \\ A^{(n)}\end{bmatrix}\,.
\end{equation}
The matrix $S_{k+1}^{-1}$ is almost diagonal. Namely, 
\[S_1^{-1}=\diag\Big(\frac1{[s_1^{(1)}]^2},\frac1{[s_1^{(2)}]^2},\ldots,\frac1{[s_1^{(n)}]^2}\Big)\,,\]
and by (\ref{nondiag}) we update its elements with
\[\frac{A^{(i)}A^{(j)}}{\tau^2}=\frac{4\Delta_k^{(i)}\Delta_k^{(j)}c_k^{(i)}c_k^{(j)}}{\tau^2[\sigma^{(i)}]^2[\sigma^{(j)}]^2}\,.\]
$c_k$ is quickly converging and $\Delta_k^{(i)}\Delta_k^{(j)}$ (for $i\ne j$) has a Bernoulli $\pm 1$ distribution with probability $\frac12$ for each $\pm 1$ outcome,
which means that non-diagonal elements remain pretty close to zero all the time if we take $\sigma^{(i)},\sigma^{(j)}$ large enough (this can always be achieved by rescaling $\theta$ if necessary). This is the reason why we used Gauss-Jordan elimination method without pivoting \cite[Section 4.2]{cia} for solving the system (\ref{parupd}) and never encountered neither the speed problems, neither the division by zero problems.

Finally, we update our belief about parameters simply with
\begin{equation}\label{bspsaupd}
\theta_{k+1}^{(i)}=\theta_k^{(i)}+b_k^{(i)}\,.
\end{equation}

\section{Experimental results}

As already announced, we now compare BSPSAS (BSPSA simple) method given by (\ref{tindup}),(\ref{sindup}) and BSPSA method given by (\ref{nondiag})--(\ref{bspsaupd}) with SPSA method described in \cite{spc} and implemented in \cite{spsa}. SPSA updates parameters by formula
\[\theta_{k+1}^{(i)}=\theta_k^{(i)}+\frac{a_k^{(i)}}{\Delta_k^{(i)}c_k^{(i)}}\cdot w_k\,, \]
where
\begin{equation}\label{spsapupd}
a_k^{(i)}=\frac{a^{(i)}}{{(A+k)}^\alpha},\quad c_k^{(i)}=\frac{c^{(i)}}{k^\gamma}\,, 
\end{equation}
and where we use $\alpha=0.602$, $\gamma=0.101$, $A=0.1N$ ($N=$ total number of iterations) -- values recommended in \cite{spc}.
Actually, in \cite{spsa} new variables are introduced by
\[R_k^{(i)}=\frac{a_k^{(i)}}{[c_k^{(i)}]^2} \]
and $R^{(i)}=R_N^{(i)}$ are used as input variables. We use this implementation, and so we discuss optimal values for $R^{(i)}$ later on.

In BSPSA(S), we have used the same formula for $c_k^{(i)}$ so far, but that is also a somewhat problematic choice, as it will be discussed later. We also use fixed value $\tau=0.6$. Let us recall that $\tau$ is the standard deviation of the approximate conditional distribution of $w_k$ (a result of the two-game match), which can be calculated to be $0.6$ for two equally strong engines with $82\%$ draw rate (which is measured to be the case in our experimental tests). So, for useful application of BSPSA(S) it remains to obtain optimal values for $s_1^{(i)}$ and $\sigma^{(i)}$, just as it remains to obtain optimal values for $R^{(i)}$ in SPSA method. Before we start, we recall the famous Elo method for measuring the relative strength of chess engines/players (see \cite{elo} for more details).
\begin{definition}\label{elo}
We say that an engine/player $E_1$ is $x$ Elo points stronger than an engine/player $E_2$ if his/her expected score (win percentage $WP$ in our earlier notation) in their match is
\[WP=\frac1{1+10^{-\frac{x}{400}}}\,. \]
We also use the notation
\[\Elo(x)=|p_2-p_1| \]
if an engine with parameter $p_2\in\R$ is $x$ Elo points stronger than an engine with parameter $p_1\in\R$.
\end{definition}

First, we did many experiments under the perfect conditions of the simulator. Namely, we assumed that all parameters have an optimal value $0$ and that Elo loss is quadratic around the optimum, which in the notation of the previous definition means $x=a(p_1^2-p_2^2)$ for some constant $a>0$. The code written in Perl \cite{wco} for the BSPSA algorithm (simulator mode is the part of the code), and the simulator results are available at \cite{bspsa}.

Our first goal was to find hyperparameters ($R^{(i)}$ for SPSA and $s_1^{(i)},\sigma^{(i)}$ for BSPSA(S)) that give the best results (on average) in the simulator mode. For SPSA we got an approximation
\begin{equation}\label{spsacon}
R^{(i)}=\frac{19362\cdot\ln\big(1+\frac{\Elo(100)}{11405}\big)}{N^{0.6} (c_N^{(i)})^{1.6}}\,, 
\end{equation}
where, as already noted above, $N$ is the total number of iterations during the execution of the algorithm. For BSPSA(S) we got much simpler formulae:
\begin{equation}\label{bspsacon}
s_1^{(i)}=\Delta\theta^{(i)},\quad \sigma^{(i)}=\Elo(100)\,,
\end{equation}
which do not depend on the number of iterations $N$, but which in turn require estimation of $\Delta\theta^{(i)}$ -- the distance of the parameter $\theta^{(i)}$ from the optimum. Furthermore, in both cases, we also need to estimate $\Elo(100)$, as defined in Definition \ref{elo}. That can be done with high precision with a relatively low number of games, and during that calculation, we can also get some approximation for $\Delta\theta^{(i)}$.

Then we used formulae (\ref{spsacon}) and (\ref{bspsacon}) to run simulations with $1,4,16$ and $64$ parameters and $200000$ iterations, where in all cases, the total distance from the optimum was $2$ Elo points at the beginning of the simulation. We discuss the choice for $c_N^{(i)}$ later in the concluding section. We repeated each run 50 times to reduce statistical errors. The obtained results are presented in the following table.

\begin{center}
 \begin{tabular}{||c || c | c | c | c ||} 
\hline
 \multirow{2.4}{*}{Method} & \multicolumn{4}{ c ||}{Elo gain mean / standard deviation}  \\ [0.3ex]
\cline{2-5}
 & 1 parameter & 4 parameters & 16 parameters & 64 parameters \\ [0.3ex]
 \hline\hline
 SPSA & 1.99944 / 0.00071 & 1.9929 / 0.0044 & 1.9048 / 0.0331 & 1.2616 / 0.1159 \\ 
 \hline
 BSPSAS & 1.99970 / 0.00042 & 1.9951 / 0.0035 & 1.9333 / 0.0256 & 1.2525 / 0.1093 \\
 \hline
 BSPSA & 1.99968 / 0.00044 & 1.9953 / 0.0033 & 1.9303 / 0.0243 & 1.2683 / 0.1054 \\
 \hline
\end{tabular}
\end{center}
We see that under perfect simulator conditions when parameters are independent and parabolically behaved around the optimum, there is basically no difference between BSPSAS given by (\ref{tindup}),(\ref{sindup}) and BSPSA given by (\ref{nondiag})--(\ref{bspsaupd}). So, for tuning a large number of parameters ($>100$), where speed could be crucial, BSPSAS seems like a great choice, but for tuning a lower number of parameters, we still prefer BSPSA. We see that BSPSA(S) comes closer to the optimum (up to more than 40\%) than SPSA, and the closer to the optimum we are, the more significant the difference is. BSPSA(S) is also more stable, with a lower standard deviation of the measurements. The difference is low for $64$ parameters, but it seems that 200000 iterations are simply not enough in that case.

Finally, we put our method to the real test. We modified $29$ parameters in the Stockfish chess engine to obtain an engine that is around $180$ Elo weaker than the original Stockfish. Then we tried to tune those parameters with 100000 iterations using both SPSA and BSPSA and using formulae (\ref{spsacon}) and (\ref{bspsacon}) to calculate input values for each tuner. Now the real chess games were played with $5+0.05$ time control ($5$ seconds for the whole game plus $50$ milliseconds increment after each move -- the maximum we could afford with our hardware). After tuning with SPSA, an engine was around $64$ Elo weaker than the original Stockfish, while tuning with BSPSA brought him to $-50$ Elo compared to the original Stockfish, which is more than $20\%$ closer to the optimum than the result obtained with SPSA. Moreover, we noticed that one parameter got stuck in the local optimum. After modifying that single parameter to its original value, we were around $-22$ Elo points for SPSA and around $-15$ Elo points for BSPSA.

\section{Conclusion}

An improvement was found over the classical SPSA method. Moreover, we believe that further improvements are possible because we used the same rule (\ref{spsapupd}) for updating $c_k^{(i)}$ which has some side effects. Namely, if $c_k^{(i)}$ values are too large and parameter behavior around an optimum is asymmetrical, we shall miss the optimum no matter how many iterations our algorithm takes. On the other hand, if $c_k^{(i)}$ values are too small, an algorithm will have slow progress and likely get stuck in the local optimum. The following picture shows an Elo dependence of one parameter that got stuck in the red region during our tuning tests.

\begin{figure}[H]
\[\includegraphics[width=0.333\textwidth]{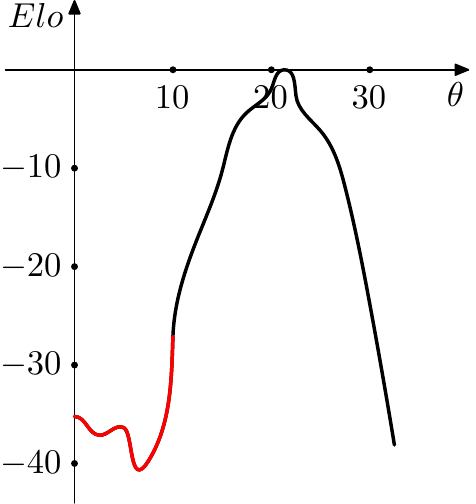}\]
\vspace{-0.9cm}
\caption*{Behaviour of one critical parameter}
\vspace{-0.4cm}
\end{figure}

Here is an idea for future work: to develop an algorithm that will change $c_k^{(i)}$ values intelligently, trying to distinguish between local and global optima. We believe that such algorithms could then be applied in many research areas, including training of neural networks and some strictly deterministic problems like finding optima of functions of several (many) variables.

\end{document}